\documentclass[12pt,reqno]{amsart}
\usepackage{amscd,amssymb,amsmath,color}
 \usepackage[all]{xy}
\numberwithin{equation}{section}
\usepackage{amsfonts}
\usepackage{amsthm}
\usepackage{color}
\usepackage{float}
\usepackage{hyperref}
\usepackage[dvipsnames]{xcolor}

\usepackage{braket}

\usepackage{todonotes}
\usepackage{blkarray}
\usepackage{pgf,tikz}
\usepackage{mathrsfs}
\usetikzlibrary{arrows}
\definecolor{cqcqcq}{rgb}{0.7529411764705882,0.7529411764705882,0.7529411764705882}
\definecolor{qqzzff}{rgb}{0.,0.6,1.}
\definecolor{wwccqq}{rgb}{0.4,0.8,0.}
\definecolor{ffqqqq}{rgb}{1.,0.,0.}

\input{prepictex.tex}
\input{pictex.tex}
\input{postpictex.tex}

\usepackage{todonotes}
\usepackage{pgf,tikz}

\title{A vector space basis of the quantum symplectic sphere}

\author{Sophie Emma Zegers} \thanks{This work was  supported by  the 
DFF-Research Project 2 on `Automorphisms and invariants of operator algebras', Nr. 7014--00145B.}
\address{Department of Mathematics and Computer Science, University of Southern Denmark, 
Campusvej 55, 5230 Odense M, Denmark} 
\email{mikkelsen@imada.sdu.dk} 

\date\today



\begin{document}

\begin{abstract}
We present a candidate of a vector space basis for the algebra $\mathcal{O}(S_q^{4n-1})$ of the quantum symplectic sphere for every $n\geq 1$. The algebra $\mathcal{O}(S_q^{4n-1})$ is defined as a certain subalgebra of the quantum symplectic group $\mathcal{O}(SP_q(2n))$. A non-trivial application of the Diamond Lemma is used to construct the vector space basis and the conjecture is supported by computer experiments for $n=1,2,...,8$. 
\end{abstract}

\theoremstyle{plain}
\newtheorem{theorem}{Theorem}[section]
\newtheorem{corollary}[theorem]{Corollary}
\newtheorem{lemma}[theorem]{Lemma}
\newtheorem{proposition}[theorem]{Proposition}
\newtheorem{conjecture}[theorem]{Conjecture}
\newtheorem{commento}[theorem]{Comment}
\newtheorem{problem}[theorem]{Problem}
\newtheorem{remarks}[theorem]{Remarks}

\theoremstyle{definition}
\newtheorem{example}[theorem]{Example}

\newtheorem{definition}[theorem]{Definition}

\newtheorem{remark}[theorem]{Remark}

\newcommand{\Nb}{{\mathbb{N}}}
\newcommand{\Rb}{{\mathbb{R}}}
\newcommand{\Tb}{{\mathbb{T}}}
\newcommand{\Zb}{{\mathbb{Z}}}
\newcommand{\Cb}{{\mathbb{C}}}

\newcommand{\Af}{\mathfrak A}
\newcommand{\Bf}{\mathfrak B}
\newcommand{\Ef}{\mathfrak E}
\newcommand{\Gf}{\mathfrak G}
\newcommand{\Hf}{\mathfrak H}
\newcommand{\Kf}{\mathfrak K}
\newcommand{\Lf}{\mathfrak L}
\newcommand{\Mf}{\mathfrak M}
\newcommand{\Rf}{\mathfrak R}

\newcommand{\x}{\mathfrak x}
\def\C{\mathbb C}
\def\N{\mathbb N}
\def\R{\mathbb R}
\def\T{\mathbb T}
\def\Z{\mathbb Z}

\def\A{{\mathcal A}}
\def\B{{\mathcal B}}
\def\D{{\mathcal D}}
\def\E{{\mathcal E}}
\def\F{{\mathcal F}}
\def\G{{\mathcal G}}
\def\H{{\mathcal H}}
\def\J{{\mathcal J}}
\def\K{{\mathcal K}}
\def\LL{{\mathcal L}}
\def\N{{\mathcal N}}
\def\M{{\mathcal M}}
\def\N{{\mathcal N}}
\def\OO{{\mathcal O}}
\def\P{{\mathcal P}}
\def\Q{{\mathcal Q}}
\def\SS{{\mathcal S}}
\def\T{{\mathcal T}}
\def\U{{\mathcal U}}
\def\W{{\mathcal W}}

\def\ext{\operatorname{Ext}}
\def\span{\operatorname{span}}
\def\clsp{\overline{\operatorname{span}}}
\def\Ad{\operatorname{Ad}}
\def\ad{\operatorname{Ad}}
\def\tr{\operatorname{tr}}
\def\id{\operatorname{id}}
\def\en{\operatorname{End}}
\def\aut{\operatorname{Aut}}
\def\out{\operatorname{Out}}
\def\coker{\operatorname{coker}}
\def\pol{\mathcal{O}}

\def\la{\langle}
\def\ra{\rangle}
\def\rh{\rightharpoonup}
\def\cl{\textcolor{blue}{$\clubsuit$}}

\def\bst{\textcolor{blue}{$\bigstar$}}

\newcommand{\Aut}[1]{\text{Aut}(#1)}
\newcommand{\Fn}{\mathcal{F}_n}
\newcommand{\On}{\mathcal{O}_n}
\newcommand{\Dn}{\mathcal{D}_n}
\newcommand{\Sn}{\mathcal{S}_n}

\newcommand{\norm}[1]{\left\lVert#1\right\rVert}
\newcommand{\inpro}[2]{\left\langle#1,#2\right\rangle}

\maketitle

\addtocounter{section}{-1}


\section{Introduction}

The base of the field of noncommutative geometry is the equivalence between the categories of commutative $C^*$-algebras and locally compact Hausdorff spaces due to Gelfand. Many researchers have been working on translating classical concepts into the noncommutative world by studying noncommutative algebras and $C^*$-algebras. Moreover, many classical spaces have been given a quantum analogue in terms of noncommutative algebras and their enveloping $C^*$-algebras. With the Gelfand duality in mind, one often thinks of quantum analogues of classical spaces as algebras of continuous functions on a non-existing virtual space.

The quantum sphere by Vaksman and Soibelman \cite{vs} is well-known in noncommutative geometry and is an object of great interest in several applications. Due to Hong and Szymański \cite{HS02}, the quantum, sphere by Vaksman and Soibelman can be described as a graph $C^*$-algebra, from which we can obtain useful informations about the $C^*$-algebra by considered the underlying directed graph (see e.g. \cite{bpis,flr}).

In \cite{LPR06}, another quantum version of the classical sphere, called the quantum symplectic $(4n-1)$-sphere, denoted $S_q^{4n-1}$, $q\in (0,1)$, is investigated as the total space of a noncommutative version of the classical Hopf bundle. The $*$-algebra of the quantum symplectic sphere, denoted  $\OO(S_q^{4n-1})$, is defined as a certain subalgebra of the quantum symplectic group $\mathcal{O}(SP_q(2n))$ \cite{FRT90, KS97} and is given by a set of generators and relations. Its enveloping $C^*$-algebra is denoted $C(S_q^{4n-1})$ and is thought of as the algebra of continuous functions on the virtual space $S_q^{4n-1}$. It follows by the defining relations that if $q=1$ then $C(S_q^{4n-1})$ is precisely the commutative algebra of continuous functions on the classical sphere $S^{4n-1}$. 

In \cite{Sa17a} the quantum symplectic sphere is investigated further under the name \textit{quantum quaternion sphere} and is denoted $H_q^{2n}$. It is shown that the quantum homogeneous space $C(SP_q(2n)/SP_q(2n-2))$ is isomorphic to $C(H_q^{2n})$. Moreover, all irreducible representations of  $C(H_q^{2n})$ are determined by using the ones of $C(SP_q(2n))$. The irreducible representations of $C(SP_q(2n))$ are classified due to a more general result by Korogodski and Soibelman \cite{KS98}. In \cite{Sa17b}, it is shown that $C(H_q^{2n})$ is isomorphic to the algebra of continuous functions on the Vaksman and Soibleman quantum $(4n-1)$-sphere.  Hence, at the $C^*$-algebraic level the quantum spheres coincide, which is not the case at the purely algebraic level. 
\\

In the present paper we construct a candidate of a vector space basis (see Conjecture \ref{vsbasis}) for the $*$-algebra $\OO(S_q^{4n-1})$ of the quantum symplectic sphere using the Diamond Lemma by Bergman \cite{b78}.  In the case where $n=2$ such a vector space basis has already been described in \cite[Section 2A]{Pa05}. The conjecture is supported by crucial  computer experiments for $n=1,2,...,8$ and if the conjecture is indeed true for a given $n$, we obtain a useful vector space basis of  $\OO(S_q^{4n-1})$ consisting of monomials.  

An application of the Diamond Lemma requires a great amount of calculations in showing that a certain class of monomials called \textit{ambiguities} are \textit{resolvable}. An explicit formula shows that the number of ambiguities increases polynomially as a function on $n$. This makes calculations by hand extremely time consuming. The case where $n=2$ has been calculated by hand by the author and a set-up to apply the Diamond Lemma has been described (see \cite[Section 3.6.2]{M21}). We will show that the set-up to apply the Diamond Lemma can be extended from $n=2$ to an arbitrary $n$. The difficulty, coursed by the drastically increase of the number of ambiguities, is overcome by writing a computer program.  For a fixed $n$, the program calculates all the ambiguities and determines whether they are resolvable or not. By applying the program, we conclude that Conjecture \ref{vsbasis} is indeed true for $n=1,2...,8$. From these data the practical running time of the program is estimated to increase polynomially as a function on $n$. For interest in a particular $n$, it is possible, by the set-up given in the present paper and the program, to determine if Conjecture \ref{VSbasis} is indeed true. Moreover, the program can be used to express a given monomial in the basis. 

The structure of the paper is as follows: In Section \ref{QSsphere} we first describe the object of interest, namely the algebra of the quantum symplective sphere. Then, in Section \ref{DiamondLemma}, we recall the Diamond Lemma, and apply it on the quantum symplectic sphere to set up a conjecture of a vector space basis in Section \ref{VSB}. Finally, in Section \ref{ambiguities} we describe the Python program which confirms the conjecture for at least $n=1,2,...,8$. The program is available at the authors webpage (\href{https://sophiemath.dk/research}{sophiemath.dk/research}). 


\section{The Quantum symplectic $(4n-1)$-sphere}\label{QSsphere}
In the present paper all algebras will be over the complex numbers. 
Let $q\in (0,1)$, the following defines the quantum symplectic $(4n-1)$-sphere for $n\in\Nb, n\geq 1$ (see \cite{LPR06, Sa17a}) which is a subalgebra of the algebra of the quantum symplectic group $\OO(SP_q(2n))$. The algebra $\OO(SP_q(2n))$ is defined by the symplectic $R$-matrix in \cite{FRT90} (see also \cite[Section 9.3]{KS97}). 
\begin{definition}
$\OO(S_q^{4n-1})$ is the $*$-algebra generated by the elements $\{x_i,y_i\}_{i=1}^n$ and their adjoints, subject to the relations: 
\begin{equation}\label{relation1}
\begin{aligned}
&x_ix_j=q^{-1}x_jx_i \ \forall i<j, \ \ y_iy_j=qy_jy_i \ \forall j>i, \ \ y_jx_i=qx_iy_j \ \forall i\neq j,  
\end{aligned}
\end{equation}
\begin{equation}\label{xiyi}
\begin{aligned}
& y_ix_i=q^{2}x_iy_i+(q^2-1)\sum_{k=1}^{i-1}q^{i-k}x_ky_k
\end{aligned}
\end{equation}
\begin{equation}\label{relation2}
\begin{aligned}
 \ \ &x_ix_i^*=x_i^*x_i+(1-q^2)\sum_{k=1}^{i-1}x_k^*x_k, \\
&y_iy_i^*=y_i^*y_i+(1-q^2)\left(q^{2(n+1-i)}x_i^*x_i+\sum_{k=1}^n x_k^*x_k+ \sum_{k=i+1}^n y_k^*y_k\right),
\end{aligned} 
\end{equation}
\begin{equation}\label{relation3}
\begin{aligned}
&x_iy_i^*=q^2y_i^*x_i, \ \ x_ix_j^*=qx_j^*x_i \ \forall i\neq j, \\[10pt] &y_iy_j^*=qy_j^*y_i-(q^2-1)q^{2n+2-i-j}x_i^*x_j \ \forall i\neq j , \ \ 
\\[10pt]
&x_iy_j^*=qy_j^*x_i \ \forall i<j, 
\ \ x_iy_j^*=qy_j^*x_i+(q^2-1)q^{i-j}y_i^*x_j \ \forall i>j. \\[10pt]
\end{aligned} 
\end{equation}
and the sphere relation 
$$
\sum_{i=1}^n(x_i^*x_i+y_i^*y_i)=1. 
$$
\end{definition}
Since the sphere relation is satisfied the generators will in any norm be restricted by 1, hence the enveloping $C^*$-algebra of $\OO(S_q^{4n-1})$ indeed exists and is denoted $C(S_q^{4n-1})$. 

Note that if $n=1$ then $\OO(S_q^3)$ is the same as $\OO(SU_{q^2}(2))$ defined by Woronowicz \cite{w1}, for which a vector space basis is already known \cite[Theorem 1.2]{w1}. 

\begin{remark} Our notation differs from the notation used in \cite{LPR06} and  \cite{Sa17a} . We can pass to the notation in \cite{LPR06} by denoting $y_i$ by $x_{2n+1-i}$ for $i=1,2,...,n$ and replacing $q$ with $q^{-1}$. To obtain the notation from \cite{Sa17a} we denote $y_i$ with 	$z_i^*$ and $x_i$ with $z_{2n+1-i}^*$ for $i=1,2,...,n$.
\end{remark}

\section{The Diamond Lemma}\label{DiamondLemma}
We will state the Diamond Lemma as presented in \cite{b78} in the case where the commutative associated ring is the complex numbers. 
Let $X$ be a set and denote by $\left<X\right>$ the free unital semigroup generated by $X$ with unit denoted $1$. Denote by $\C\left<X\right>$ the free algebra over $\C$ generated by $X$. 
We define a \textit{reduction system} $S$ as follows. An element $\sigma\in S$ is a tuple which we denote by $(w_{\sigma},f_{\sigma})$ with $w_{\sigma}\in \left<X\right>$ and $f_{\sigma}\in \C\left<X\right>$.  

Let $\sigma\in S$ and $a,b\in \left<X\right>$. A \textit{reduction} is an endomorphism $r_{a\sigma b}: \C\left<X\right>\to \C\left<X\right>$ which fixes all elements of $\left<X\right>$ except $aw_\sigma b$ which is mapped to $af_{\sigma} b$. A reduction $r_{a\sigma b}$ is said to act trivially on $x\in \C\left<X\right>$ if the coefficient $aw_{\sigma}b$ in $x$ is zero. 

\begin{figure}[H]
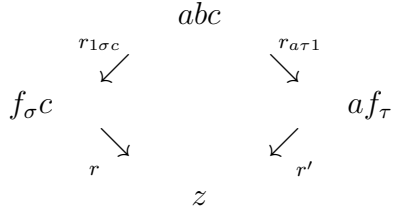

$$
\begin{matrix}
& & abc & &  \\
&\overset{r_{1\sigma c}\ \ \ }{\swarrow} & & \overset{ \ \ \ r_{a\tau 1}}{\searrow} &  \\
f_{\sigma}c & & &  & af_{\tau} \\
& \underset{r \ \ \ \ }{\searrow} &  & \underset{ \ \ \ \ r'}{\swarrow} &  \\
& & z & & 
\end{matrix}
$$
\caption{Overlap ambiguity}
\label{Overlap ambiguities}
\end{figure}

An \textit{(overlap) ambiguity} of $S$ is a $5$-tuple $(\sigma, \tau, a,b,c)$ where $\sigma, \tau\in S$ and $a,b,c\in \left<X\right>$ with $a,b,c\neq 1$ such that
$w_{\sigma}=ab$ and $w_{\tau}=bc$. 
An overlap ambiguity is \textit{resolvable} if $f_{\sigma}c$ and $af_{\tau}$ can be reduced to a common expression using compositions of reduction maps. This is illustrated in Figure \ref{Overlap ambiguities}. The arrows in the top indicates the application of the reduction $\sigma$ and $\tau$. The arrows denoted by $r$ and $r'$ indicates the composition of reduction maps we have to apply to get a common expression denoted $z\in \C \left<X\right>$. 
We wish to have a partial order on $ \left<X\right>$ such that if $b<b'$ then $abc<ab'c$ for $a,b,b',c\in  \left<X\right>$. The partial order is \textit{compatible} with the reduction system $S$, if for any $\sigma=(w_\sigma,f_{\sigma})$ in $S$, with 
$$
f_{\sigma}=\sum \alpha_j a_j, \ \ a_j\in  \left<X\right>, \alpha_j\in \C,
$$
we have $a_j<w_\sigma$ for all $j$. 
A partial order has the \textit{descending chain condition} if there are no infinite descending chains.  \\

An element $x\in \C\left<X\right>$ is called \textit{irreducible} if every reduction is trivial on $x$, in that case $x$ involves no monomials of the form $aw_{\sigma}b$ with $a,b\in  \left<X\right>$ and $\sigma\in S$.  The set of irreducible elements of $\C \left<X\right>$ is denoted $\C \left<X\right>_{irr}$. 

Let $I$ be the 2-sided ideal of $\mathbb{C}\left<X\right>$ generated by $w_{\sigma}-f_{\sigma}$ for all $\sigma\in S$. We can then state the Diamond Lemma. 
\begin{theorem}[{\cite[Theorem 1.2]{b78}}] \label{diamondlemma}

Let $\leq$ be a partial order on $ \left<X\right>$ with the descending chain condition. If $S$ is a reduction system on $\C \left<X\right>$ compatible with $\leq$ such that all ambiguities are resolvable, then a vector space basis for $\C \left<X\right>/I$ is given by $\C \left<X\right>_{irr}$. 
\end{theorem}

\section{A vector space basis for $\OO(S_q^{4n-1})$}\label{VSB}
We will now describe how to set up a proof of the following conjecture, using the Diamond Lemma. 
\begin{conjecture}\label{VSbasis}
The following elements 
\begin{equation}\label{vsbasis}
\begin{aligned} 
&{y_1^*}^{m_1}\cdots {y_n^*}^{m_n}{x_1^*}^{k_1}\cdots {x_n^*}^{k_n}x_{n-1}^{s_{n-1}}\cdots x_1^{s_1}y_n^{t_n}\cdots y_1^{t_1}, \\
&
{y_1^*}^{m_1}\cdots {y_n^*}^{m_n}{x_1^*}^{k_1}\cdots {x_{n-1}^*}^{k_{n-1}}{x_n}^{s_n}\cdots x_1^{s_1}y_n^{t_n}\cdots y_1^{t_1},
\end{aligned} 
\end{equation}
with $m_j,k_j,s_j,t_j\geq 0, y_j^0={y_j^*}^0=x_j^0={x_j^*}^0=1$, form a vector space basis for $\OO(S_q^{4n-1})$. 
\end{conjecture}
We can set up the initiation of a proof for a general $n$, but when we have to show that all ambiguities are resolvable we will no longer consider a general $n$. Instead, in Section \ref{ambiguities}, we perform computer experiments for $n=1,2,3,...,8$ supporting Conjecture \ref{vsbasis}. 
\\

In the following we describe a reduction system and a partial order which has the required properties to apply the Diamond Lemma. 

Let $X=\{x_1,...,x_n, x_1^*,...,x_n^*,y_1,...,y_n, y_1^*,...,y_n^*\}$.  
We define a reduction system, denoted $S$, by the relations in \eqref{relation1},\eqref{xiyi},\eqref{relation2},\eqref{relation3} and the following rewritten form of the sphere relation: 
$$
x_n^*x_n=1-\sum_{i=1}^{n-1}x_i^*x_i-\sum_{i=1}^n y_i^*y_i. 
$$
An element in $S$ is then a tuple which contains the left and the right hand side of these relations e.g. $(x_1y_3,q^{-1}y_3x_1)\in S$. Let $I$ be the 2-sided ideal of $\mathbb{C}\left<X\right>$ generated by $w_{\sigma}-f_{\sigma}$ for all $\sigma\in S$. Then $\OO(S_q^{4n-1})$ equals $\C\left<X\right>/I$. We can then obtain a vector space basis for $\OO(S_q^{4n-1})$ if the assumptions of Theorem \ref{diamondlemma} are satisfied. 

We begin showing that there 
exists a partial order on $\left<X\right>$ which is compatible with the reduction system. We can define a relation on $\left<X\right>$ as follows: 
\begin{definition}
For $u,v\in \left<X\right>$, we let $u<v$ if we can write $v$ as a sum $\sum_{j} \alpha_ja_j$ with $a_j\in \left<X\right>, \alpha_j\in \C$ and $u=a_j$ for some $j$, by using a finite number of the reductions in $S$ . 
\end{definition}
By construction the relation is compatible with the reduction system. Intuitively, $u<v$ if $u$ requires less reductions to be of the form (\ref{vsbasis}) than $v$. Note, it follows immediately that if $b<b'$ then $abc<ab'c$ for $a,b,b',c\in  \left<X\right>$. 

The relation $<$ is by definition reflexive and transitive. To show it is antisymmetric we will prove that $u<v$ and $v<u$ cannot happen for $u,v\in \left<X\right>$. This will be done by associating some numbers to each monomial. 

Let $u=a_{1}a_{2}\cdots a_{k}$ be a monomial in $\left<X\right>$ i.e. $a_i\in\{x_j,x_j^*,y_j,y_j^*| j=1,...,n\}$. We say that $(a_{s}, a_{t})$ is an inverted pair in $u$ if $a_{s}$ does not precede $a_{t}$ in the list 
\begin{equation}\label{order}
y_1^*y_2^*\cdots y_n^*x_1^*x_2^*\cdots x_n^* x_n\cdots x_2 x_1y_n\cdots y_2 y_1.
\end{equation}
That is, if their order as factors of $u$ is opposite that of their order in (\ref{order}). 
For each monomial $u\in \left<X\right>$ we associate numbers as follows: 
\begin{itemize}
\item[(1)] $N_1(u):=$ the number of inverted pairs of the form $(x_i,x_j^*),(x_i,y_j^*)$ and $(y_j,y_i^*)$ for $i,j=1,...,n$ in $u$. 

\item[(2)] $N_2(u):=$ minimum of the number of $x_n$ and the number of $x_n^*$ in $u$. 
\item[(3)] For $i=1,...,n$ let:  $N_{3,i}(u):=$ the number of inverted pairs of the form $(y_i, x_i)$ and $(x_i^*,y_i^*)$ in $u$. 
\item[(4)]  $N_{4}(u):=$ the number of inverted pairs of the form  $(y_j,x_i)$ and $(x_i^*,y_j^*)$ for $i,j=1,...,n$, $i\neq j$ in $u$. 
\item[(5)] $N_{5}(u):=$ the number of inverted pairs of the form $(y_i, y_j)$ and $(y_j^*, y_i^*)$ for $i=1,...,n-1$, $j=i+1,...,n$ in $u$.  
\item[(6)]  $N_{6}(u):=$ the number of inverted pairs of the form $(x_i, x_j)$ and $(x_j^*, x_i^*)$ for $i=1,...,n-1$, $j=i+1,...,n$ in $u$. 
\end{itemize}

\begin{example}
Let $n=3$. Consider the monomial 
$$
u=x_3y_3x_3x_1^*y_1y_1x_3^*x_1
$$
then 
$$
N_1(u)=8, \ N_2(u)=1, \ N_{3,1}(u)=2, \ N_{3,2}=0, \ N_{3,3}=1, \ N_4(u)=4, \ \ N_5(u)=0, \ \ N_6(u)=2. 
$$
\end{example} 

The inverted pairs counted by $N_k$, $k=1,4,5,6$ and $N_{3,i}$, $i=1,...,n$ are said to be of type $(k)$ and $(3,i)$ respectively. For $(w_{\sigma},f_{\sigma})\in S$ we say that $w_{\sigma}=ab$, $a,b\in\{x_i,y_i,x_i^*,y_i^*| \ i=1,...,n\}$ is of type $(k)$ or $(3,i)$ if $(a,b)$ is of type $(k)$ or $(3,i)$ respectively. 
 
Let $\sigma=(w_{\sigma},f_{\sigma})\in S$, then we have two different cases: 
\begin{itemize}
\item[] \textit{Case 1}: $f_{\sigma}$ is a scalar multiple of the given pair written in the switched order. 
\item[] \textit{Case 2}: $f_{\sigma}$ is a linear combination of at least two monomials. 
\end{itemize}
Let $v\in \left<X\right>$, denote by $u<v$ one of the polynomials occurring after a reduction.  Let $\sigma\in S$ with  $w_{\sigma}$ of type $(1)$. Both in case 1 and 2 we have $N_1(u)=N_1(v)-1$. 
If we are in case 1, we furthermore have  $N_k(u)=N_k(v)$, $k=2,4,5,6$ and $N_{3,i}(u)=N_{3,i}(v), i=1,...,n$. In case 2, it might happen that these numbers will increase. 

If we continue using reductions of type (1), we will at some point arrive at $N_1=0$. The monomials will then be such that all $x_i^*,y_j^*, i,j=1,...,n$ are on the left hand side and all the $x_i,y_j, i,j=1,...,n$ will be on the right hand side. 

Consider now the number $N_2(v)$. When using the reduction
$$\sigma=\left(x_n^*x_n,1-\sum_{i=1}^{n-1}x_i^*x_i-\sum_{i=1}^n y_i^*y_i\right), 
$$
coming from the sphere relation, 
we get $N_2(u)=N_2(v)-1$ and $N_1(u)=N_1(v)$. In case 1 we furthermore have that all the other numbers will be the same for $u$ and $v$. If we apply reductions of this form until $N_2=0$ we obtain a monomial which at most contains a single $x_n$ or a single $x_n^*$. 
\\

Let $\sigma\in S$ be a reduction with $w_{\sigma}$ of type $(3,m)$ i.e. coming from the relation in  (\ref{xiyi}). Then 
$$N_{3,m}(u)=N_{3,m}(v)-1, \ \ N_{3,i}(u)=N_{3,i}(v), i>m
$$
since if $w_{\sigma}$ is of type $(3,m)$ we will not create any new terms containing any $x_i, y_i$ or $x_i^*,y_i^*$ with $i>m$. Note that it might happen that $N_{3,i}(u)>N_{3,i}(v)$ if $m>i$.  

Furthermore, we have 
$$
N_1(u)=N_1(v) \ \text{and} \ N_2(u)=N_2(v). \ \ 
$$

Finally, we let $\sigma\in S$ be such that $w_{\sigma}$ is of type $(i)$ for $i=4,5,6$. Then 
$$
N_{i}(u)=N_{i}(v)-1, \ N_{l}(u)=N_{l}(v), l\neq i
$$
and 
$$
N_1(u)=N_1(v), \ N_2(u)=N_2(v), \ N_{3,m}(v)=N_{3,m}(u), m=1,...,n. 
$$
If we reduce $N_4$ to zero we will have ordered the monomials such that we first have all the $y_i^*$'s then the $x_j^*$'s afterwards $x_s$'s and in the end the $y_t$'s. When we use reductions to decrease $N_5$ to zero we order all the $y_i$'s and $y_i^*$'s into the order in (\ref{order}). Similar for $x_i, x_j^*$ when reductions are applied to decrease $N_6$ to zero.  

This concludes that $u<v$ and $v<u$ can never occur for any pair $u,v\in \left<X\right>$, since it will contradict the properties of the numbers $N_k$ and $N_{3,i}$. The relation is then antisymmetric and hence a partial order. From the above we also get that the reduction process will terminate at some point. Hence the descending chain condition is satisfied. 
\\

If we can show that all the ambiguities are resolvable, we have by the Diamond Lemma that the elements in (\ref{vsbasis}) is a vector  spaces for $\OO(S_q^{4n-1})$. In that case it follows that $x_1,...,x_{n-1}$ are non-zero in $\OO(S_q^{4n-1})$. 
\\

The relation can also be given an alternative description as follows: 
$u<v$ if and only if one of the following conditions are satisfied: 
\begin{enumerate}
\item There exists a $k$ such that $N_k(u)<N_k(v)$ and for all $j<k$ we have  $N_j(u)=N_j(v)$. Furthermore, if $k=4,5,6 $ we must have $N_{3,i}(u)=N_{3,i}(v)$ for all $i$. 
\item There exists a $m$ such that $ N_{3,m}(u)<N_{3,m}(v)$, $N_{3,i}(u)=N_{3,i}(v)$ for $i>m$ and $N_k(u)=N_k(v)$ for $k=1,2$. 
\end{enumerate}
This gives us an efficient algorithm for reducing a monomial into a linear combination of elements of the form (\ref{vsbasis}). Applying suitable relations we reduce the number of inverted pairs of type (1). When there are no such inverted pairs (i.e. $N_1$ is $0$) we use reductions until $N_2$ is $0$ for all the monomials. Then we reduce the number of inverted pairs of type $(3,n)$, afterwards $(3,n-1)$ and so on, until we have reduced the number of inverted pairs of of type $(3,1)$. 
In the end we reduce the number of inverted pairs of type $(4), (5)$ and $(6)$ respectively. Let $n=3$, below is an example of how to use the algorithm. A different color indicates that we have moved on to another step in the reduction.
$$
\begin{aligned}
y_2\textcolor{red}{y_3x_1^*}x_3&= q\textcolor{red}{y_2x_1^*}y_3x_3 = q^2x_1^*y_2\textcolor{blue}{y_3x_3}  \\
&= q^4x_1^*y_2x_3y_3+q^4(q^2-1)x_1^*{y_2x_1}y_1 
+ q^3(q^2-1)x_1^*\textcolor{orange}{y_2x_2}y_2 \\
&= q^4x_1^*\textcolor{LimeGreen}{y_2x_3}y_3+q^4(q^2-1)x_1^*\textcolor{LimeGreen}{y_2x_1}y_1 
+ q^5(q^2-1)x_1^*x_2y_2y_2+q^4(q^2-1)^2x_1^*x_1y_1y_2 \\
&= q^5x_1^*x_3\textcolor{cyan}{y_2y_3}+q^5(q^2-1)x_1^*x_1{y_2y_1}+ q^5(q^2-1)x_1^*x_2y_2y_2+q^4(q^2-1)^2x_1^*x_1\textcolor{cyan}{y_1y_2} \\
&= q^6x_1^*x_3y_3y_2+q^5(q^2-1)x_1^*x_1{y_2y_1}+ q^5(q^2-1)x_1^*x_2y_2y_2+q^5(q^2-1)^2x_1^*x_1{y_2y_1} \\
&= q^6x_1^*x_3y_3y_2+(q^9-q^7)x_1^*x_1{y_2y_1}+ q^5(q^2-1)x_1^*x_2y_2y_2
\end{aligned}
$$


\section{Ambiguities of the reduction system}\label{ambiguities}
When the number $n$ grows larger we obtain extremely many ambiguities, making it hard to calculate by hand.  When it comes to determine whether all the ambiguities are resolvable or not, we will apply computer calculations. Another advantage using computer calculations is the precision of the calculations. It is not hard to imagine that one could make mistakes when reducing monomials by hand. 
\\

The ambiguities can be found as follows. Consider the sequence 
$$
y_1^*y_2^*\cdots y_n^*x_1^*x_2^*\cdots x_n^* x_n\cdots x_2 x_1y_n\cdots y_2 y_1.
$$
If $n=2$, take first $x_1^*$, then pick an element to the left of $x_1^*$ in the above sequence, say $y_2^*$. Pick then an element to the left of $y_2^*$, say $y_1^*$. Then we have an ambiguity $x_1^*y_2^*y_1^*$. Take now $x_2^*$, pick an element to the left, and then another element to the left of the just chosen element. Continue like this for all the elements in the above sequence. For $n> 2$, we start with $y_3^*$ and follow the same procedure as before. Furthermore, we must include the following ambiguities: 
$$
\begin{aligned}
&ax_n^*x_n, \ a\in\{x_2,...,x_1,y_n,...y_1\} , \\
&x_n^*x_nb, \ b\in\{y_1^*,...,y_n^*,x_1^*,...,x_n^*\}.  
\end{aligned}
$$
In all we get the following number of ambiguities for a specific $n$: 
$$
\begin{aligned}
\left(\sum_{i=1}^{4n-2}\sum_{k=1}^i k \right)+4n
&=\frac{1}{2}\sum_{i=1}^{4n-2}i(i+1)+4n \\
&=\frac{1}{12}(4n-2)(4n-1)(2(4n-2)+1)+\frac{1}{4}(4n-2)(4n-1)+4n
\\ 
&=\frac{1}{12}(4n-2)(4n-1)8n+4n =\frac{8}{3}(4n^3-3n^2+2n). 
\end{aligned}
$$

A computer program has been written by the author in Python Version 3.9.5 using the library SymPy, it is available at the authors webpage \href{https://sophiemath.dk/research}{sophiemath.dk/research}. The program contains the following parts: 
\begin{itemize}
\item[(1)] For a given $n$ it lists all the generators and  relations in the reduction systems. 
\item[(2)] An algorithm to reduce a given monomial until it can not be reduced any more by the reduction system.  
\item[(3)] All the ambiguities are determined and the algorithm in (2) is applied to check if the ambiguities are resolvable or not. 
\end{itemize}
Note that the monomials in the process are reduced until all the involved terms are irreducible. This is needed for the program to determine if the ambiguities are resolvable. E.g. consider the ambiguity $x_1^*y_2^*y_1$, then $x_1^*y_2^*$ and $y_2^*y_1$ are reduced until they contain only irreducible parts, afterwards they are compared. 
Note that it is not needed in theory, since it is enough to reduce until we get a common expression, which could contain reducible parts. 

For $n=1,2,3,....,8$ we obtain, by the computer program, that all the ambiguities are resolvable. Hence Conjecture \ref{vsbasis} is indeed true for $n=1,2,3,...,8$. In Table \ref{runningtime} the number of ambiguities and the practical running time\footnote{The running times has been collected by running the program on a standard laptop with 3.3 GHz Intel Core i7, dual core.}  of the program is given.

\begin{table}[H]
\resizebox{\columnwidth}{!}{
\begin{tabular}{ |c|c c c c c c c c }
$n$ & 1 & 2 & 3 & 4 & 5 & 6 & 7 & 8 \\ 
& & & & & & & \\
\hline
Number of \\
ambiguities & 8 & 64 & 232 & 576 & 1160 & 2048 & 3304 & 4992  \\
& & & & & & & \\
Running time & 4.80 & 82.29 &  509.76 & 1917.51  & 5417.61 &  12944.70 & 27948.76 & 51751.18 \\
 in seconds &  & (1.37 minutes) & (8.50 minutes) & (32.00 minutes) & (1.50 hour) & (3.60 hours) & (7.76 hours) & (14.38 hours)
\end{tabular}}
\vspace{0.3cm}
\caption{Number of ambiguities and running time in seconds.}
\label{runningtime}
\end{table}
A regression on the above data has been made. This indicates that the growth of the practical running time is likely to be polynomial. 
\\

We have then obtained computer experiments which supports Conjecture \ref{vsbasis}. For interest in the vector space basis of $\OO(S_q^{4n-1})$ for a specific $n$, this is possible to determine using the computer program, taking the running time into account. Furthermore, as a part of the program, we can express a given monomial in the basis using computer calculations. 

\begin{remark}
The algorithm in program, which checks whether the ambiguities are resolvable or not, takes as its input the generators and relations of the quantum symplectic sphere. Hence, if one is interested in a similar problem for another algebra, then it might be possible to type the new set of generators and relations and then apply the program. 
\end{remark}


\subsection*{Acknowledgment.} I would like to thank Nikolai Nøjgaard for a very helpful guidance in the process of writing the computer program.  I also gratefully acknowledge helpful comments and suggestions from Wojciech Szymański and Thomas Gotfredsen.

\end{document}